\documentclass[12pt,oneside]{article}

\usepackage{amssymb,latexsym,start2e,amsthm,amsmath}
\usepackage{pstricks,pst-node,pst-tree}
\usepackage{graphicx, pstcol, pst-plot}
\usepackage{subfigure}
\usepackage{color}
\definecolor{darkblue}{rgb}{0, 0, .4}
\definecolor{externblue}{rgb}{0, 0, .7}
\usepackage[bookmarks]{hyperref}
\hypersetup{
	colorlinks=true,
	linkcolor=darkblue,
	anchorcolor=darkblue,
	citecolor=darkblue,
	pagecolor=darkblue,
	urlcolor=externblue,
	pdftitle={Rationality of the M\"obius function \\ of a
	  composition poset
	},
	pdfsubject={Combinatorics
	},
	pdfauthor={Anders Bj\"orner and Bruce E. Sagan
	},
	pdfkeywords={automaton, composition, generating function,
	  hypergeometric series, monoid, rationality, subword order
	}
}

\newcommand{\Supp}{\operatorname{Supp}}
\newcommand{\sspn}[1]{\langle\langle{#1}\rangle\rangle}

\begin{document}

\title{Rationality of the M\"obius function \\ of a
	  composition poset
}

\author{Anders Bj\"orner\\[-5pt]
\small Department of Mathematics\\[-5pt]
\small Kungliga Tekniska H\"ogskolan\\[-5pt]
\small S-100 44 Stockholm, SWEDEN\\[-5pt]
\small
\href{mailto:bjorner@math.kth.se}{\texttt{bjorner@math.kth.se}}\\[6pt]
Bruce E. Sagan
\thanks{This work was partially done while the
author was on leave at the Institut Mittag-Leffler and at DIMACS}\\[-5pt]
\small Department of Mathematics\\[-5pt]
\small Michigan State University\\[-5pt]
\small East Lansing, MI 48824-1027, USA\\[-5pt]
\small
\href{mailto:sagan@math.msu.edu}{\texttt{sagan@math.msu.edu}}
}

\date{\today \\[6pt]
	\begin{flushleft}
	\end{flushleft}
            }

\maketitle

\begin{abstract}
We consider the zeta and M\"obius functions of a partial order on
integer compositions first studied by Bergeron,
Bousquet-M\'elou, and Dulucq.   The M\"obius function of this poset
was determined by Sagan and Vatter.  We prove rationality of various
formal power series in noncommuting variables whose coefficients 
are evaluations of the zeta function $\zeta$ and 
the M\"obius function $\mu$.  The proofs are either directly
from the definitions or by constructing finite-state
automata.  

We also obtain explicit expressions for generating
functions obtained by specializing the variables to commutative ones.
We reprove Sagan and Vatter's formula for $\mu$ using this machinery.
These results are closely related to those of Bj\"orner and Reutenauer
about subword order, and we discuss a common generalization.
\end{abstract}

\section{Introduction}

Let $A$ be an arbitrary set and consider the {\it free monoid}, $A^*$,
of all words over $A$:
$$
A^*=\{w=w(1)w(2)\ldots w(n)\ |\ \mbox{$n\ge0$ and $w(i)\in A$ for all
$i$}\}.
$$
We let $\ell(w)$ denote the length (number of elements) of $w$.

If $\bbP$ is the positive integers, then $\bbP^*$ is just the set of
integer compositions (ordered partitions).  We put a partial order on
$\bbP^*$ by saying that $u\le w$ if and only if $w$ contains a subword 
$w(i_1)w(i_2)\ldots w(i_l)$ where $l=\ell(u)$ and
$$
\mbox{$u(j)\leq w(i_j)$ for $1\leq j\leq l$.}
$$
To illustrate, $\bf 334\le 34261$  as can be seen by considering the
subword $\bf 346$.   Note that integers will be typeset in boldface
when considered as elements of $\bbP^*$.
Bergeron, Bousquet-M\'elou and Dulucq~\cite{bbd:spc} initiated the
study of $\bbP^*$ by counting its saturated lower chains.  This work
was carried on by Snellman~\cite{sne:scc,sne:spa} who also considered
saturated chains in two other partial orders on $\bbP^*$.
One of these posets was originally defined by Bj\"orner
and Stanley~\cite{bs:ayl} who showed that it has analogues of many of the
properties of Young's lattice.  Sagan and Vatter~\cite{sv:mfc}
determined the M\"obius function of the poset we are considering.
Here we will use generating functions over monoids to give more
information about the M\"obius and zeta functions of $\bbP^*$ as well
as rederiving the theorem of Sagan and Vatter using this machinery.

There is a strong connection between this order on $\bbP^*$ and
subword order.  Considering $A$ to be arbitrary, we define {\it
subword order\/} on $A^*$ by letting $u\le w$ if and only if 
there is a subword $w(i_1)w(i_2)\ldots w(i_l)$ of length $l=\ell(u)$ with
$$
\mbox{$u(j)=w(i_j)$ for $1\leq j\leq l$.}
$$
For example, $abba\le ababbbaa$ since $w(1)w(4)w(6)w(8)=abba$.
Context will make it clear whether ``$\le$'' refers to subword order
or composition order.
Bj\"orner~\cite{bjo:mfs} was the first to give a complete
characterization of the M\"obius function for subword order.
See~\cite{sv:mfc} for a history of this problem.  In particular,
Bj\"orner and Reutenauer~\cite{br:rmf}  showed that the M\"obius and
zeta functions have rational generating functions and were able to
reprove the formula for $\mu$ using these ideas.  

The rest of this paper is structured as follows.  In the next section
we provide the necessary definitions to state Bj\"orner's formula
for $\mu$ in $A^*$ as well as Sagan and Vatter's result in $\bbP^*$,
see Theorems~\ref{muA*} and~\ref{muP*}, respectively.
In Section~\ref{r} we prove the rationality of  monoid
generating functions for $\mu$ and $\ze$ on the subposet $\{\1,\2,\ldots,\bmn\}^*$
of $\bbP^*$.  Our demonstrations are either based directly on the
definitions or  use finite-state automata.  By specializing the
variables, we obtain explicit formulas 
for related generating functions in Section~\ref{gfc}.  
Surprisingly, results about hypergeometric series are needed to do
some of the computations.  The
next section is devoted to another proof of the formula
for $\mu$ in $\bbP^*$ using the generating function approach. 
Sagan and Vatter showed that both Theorems~\ref{muA*} and~\ref{muP*}
are special cases of a more general result about certain partial orders which
they called generalized subword orders (and which have been studied in
the context of well-quasi-ordering, see Kruskal~\cite{kru:twq}).  In Section~\ref{gso}, we
indicate which of our results can be proved in this level of generality.
We end with a section of comments and open problems.

\section{Subword and composition order}
\label{sco}

We will first present the formula for the M\"obius function of
$A^*$ in a way that will help motivate our definitions when we get to
$\bbP^*$.  We will not define the M\"obius function itself, but that
background can be found in the text of 
Stanley~\cite[\S  3.6--3.7]{sta:ec1}.

We begin by giving an equivalent formulation for subword order which
will be useful when we get to $\mu$.  Suppose we have a special symbol
$0$ with $0\not\in A$.  Then the {\it support\/} of a word
$\eta=\eta(1)\eta(2)\ldots\eta(n)\in (A\cup 0)^*$  is
$$
\Supp \eta =\{i\ |\ \eta(i)\neq 0\}.
$$
An {\it expansion\/} of $u\in A^*$ is a word $\eta_u\in(A\cup 0)^*$
such that the restriction $\eta_u$ to its support is $u$. 
Taking $u=abba$ as before, then one possible expansion is
$\eta_u=a00b0b0a$.  An {\it
embedding\/} of $u$ into $w$ is an expansion $\eta_u$ of $u$ 
having length $\ell(w)$ and satisfying
$$
\mbox{$\eta_u(i)=w(i)$ for all $i\in\Supp\eta_u$.}
$$
Clearly $u\le w$ in subword order if and only if there is an embedding of
$u$ into $w$.  In fact, the example $\eta_u$ above is 
the embedding which corresponds to the subword of
$w=ababbbaa$ given in the previous section.

The M\"obius function of subword order counts a particular type of
embedding.  Suppose $a\in A$.  A {\it run\/} of $a$'s in $w$ is a
maximal interval of indices $[r,t]$ such that
$$
w(r)=w(r+1)=\ldots=w(t)=a.
$$
Continuing with our example, $w=ababbbaa$ has runs
$[1,1]$, $[2,2]$, $[3,3]$, $[4,6]$, and $[7,8]$.
An embedding $\eta_u$ into $w$ is {\it normal\/} if, for every $a\in
A$ and every run $[r,t]$ of $a$'s, we have
$$
(r,t] \sbe \Supp \eta_u
$$
for the half-open interval $(r,t]$.  In our running
example, this means that the $b$'s in positions 5 and 6
as well as the $a$ in position 8
must be in any normal embedding.  (Runs of one element impose no
restriction since if $r=t$ then $(r,t]=\emp$.)
So in this case there are exactly two normal 
embeddings $\eta_u$ into $w$, namely
$$
\mbox{$\eta_u = a000bb0a$ and $00a0bb0a$.}
$$

Let $\binom{w}{u}_n$ denote the number
of normal embeddings of $u$ into $w$.
\bth[Bj\"orner~\cite{bjo:mfs}]
\label{muA*}
If $u,w\in A^*$ then\\
\eqed{\mu(u,w)=(-1)^{|w|-|u|}\dil{w\choose u}_n.}
\eth
Putting everything together in our example, we obtain
$$
\mu(abba,ababbbaa)=(-1)^{8-4} \cdot 2 = 2.
$$

In $\bbP^*$, the definitions of support and expansion are the same as
in $A^*$.  However, the definition of embedding must be changed
to reflect the different partial order.
In this case,  define an {\it embedding\/} of $u$ into
$w$ as an expansion $\eta_u$ such that $\ell(\eta_u)=\ell(w)$ and
$$
\mbox{$\eta_u(i)\le w(i)$ for $1\le i\le \ell(w)$.}
$$
As before, $u\le w$ in $\bbP^*$ if and only if there exists
an embedding of $u$ into $w$.

Of particular interest to us will be the rightmost embedding.
Suppose $u\le w$.  The {\it rightmost\/} embedding $\rho_u$ into $w$
is the one such that for any other embedding $\eta_u$ into $w$ we have
$\Supp(\rho_u)\ge \Supp(\eta_u)$.  (If $S=\{i_1<\cdots<i_m\}$ and
$S'=\{i_1'<\cdots<i_m'\}$ then $S\ge S'$  means
$i_j\ge i_j'$ for $1\le j \le m$.)  

The definition of a run is again the same in $\bbP^*$ as it was in
$A^*$.  So we call an embedding $\eta_u$ into $w$ {\it normal\/} if it
satisfies the following two criteria.
\begin{enumerate}
\item\label{normal1} For $1\le i\le \ell(w)$, we have $\eta_u(i)=w(i)$,
$w(i)-1$, or $0$.\vs{3pt}
\item\label{normal2} For all $k\ge 1$ and every run $[r,t]$ of $\bmk$'s
in $w$, we have\vs{3pt}
\begin{enumerate}
\item $(r,t]\sbe \Supp \eta_u$ if  $k=1$,\vs{3pt}
\item $r\in\Supp \eta_u$ if $k\ge 2$.
\end{enumerate}
\end{enumerate}
Note that in $\bbP^*$ a normal embedding can have three possible values at each
position instead of the two permitted in $A^*$.  Also note that the
run condition for ones is the same 
as in $A^*$, while that condition for integers greater than one is
complementary.  For example, if $u=\bf 21113$ and $w= \bf 2211133$,  
then there are two normal embeddings, namely
$\eta_u= \bf 2101130$ and $\bf 2011130$.  Also,
$\bf 2001113$  and $\bf 0211130$ are not normal since
they violate conditions~\eqref{normal1} and~\eqref{normal2},
respectively.

Another difference between $A^*$ and $\bbP^*$ is that in the former
the sign of an embedding only depends on the length difference, while
in the latter it depends on the embedding itself.
If $\eta_u$ into $w$ is normal then define its {\it defect} to be
$$
d(\eta_u)=\#\{i\ |\ \eta_u(i)=w(i)-1\}.
$$

The formula for the M\"obius function of $\bbP^*$ is as
follows.
\bth[Sagan and Vatter~\cite{sv:mfc}]
\label{muP*}
If $u,w\in\bbP^*$ then
$$
\mu(u,w)=\sum_{\eta_u} (-1)^{d(\eta_u)}
$$
where the sum is over all normal embeddings $\eta_u$ into $w$.\qqed
\eth
Finishing off the example of the previous paragraph, 
$$
\mu({\bf 21113},{\bf 2211133}) = (-1)^2+(-1)^0=2.
$$
Although this example does not show it, it is possible to have
cancellation among the terms in the sum for $\mu$.

\section{Rationality}
\label{r}

Let $\ep$ denote the empty word in $A^*$.  For this section and the
next one we will assume that $A$ is a finite set.  Let $\bbZ\sspn{A}$
be the algebra of formal power series in the noncommuting variables
$A$ with integer coefficients.  So every $f\in\bbZ\sspn{A}$ has the
form 
$$
f=\sum_w c_w\ w
$$
where $w\in A^*$ and $c_w\in\bbZ$.  If $f$ has no constant term, i.e.,
$c_\ep=0$, then define
\beq
\label{*}
f^*=\ep+f+f^2+f^3+\cdots = (\ep-f)^{-1}.
\eeq
(One needs the restriction on $f$ to make sure that the sum is well
defined as a formal power series.)  We say $f$ is {\it rational\/} if
it can constructed from a finite set of monomials using a finite
number of applications of the algebra operations and the star
operation.  For more information about rational series, see the books
of Eilenberg~\cite{eil:almB} or Berstel and Reutenauer~\cite{br:rsl}.  
We will show in this section that various series related to the
M\"obius and zeta functions are rational.  

It will be convenient to define $[n]=[1,n]$.  We will also use such
interval notations with elements of $\bbP^*$ in the obvious way.  So,
for example, 
$$
[\bmk,\bmn]=\{\bmk,\bmk\bmpl\1,\ldots,\bmn\}.
$$
Consider $[\bmn]^*$ as a
subposet of $\bbP^*$.  Given $u\in[\bmn]^*$, we have the associated
formal series
\beq
\label{Zu}
Z(u)=\sum_{w\ge u} w =\sum_{w} \zeta(u,w) w
\eeq
where $\zeta$ is the zeta function of $[\bmn]^*$.  We also wish to
consider
\beq
\label{Mu}
M(u)=\sum_{w\ge u} \left( \sum_{\eta_u} (-1)^{d(\eta_u)} \right)w
\eeq
where the inner sum is over all normal embeddings $\eta_u$ into $w$.
Note that if we assume Theorem~\ref{muP*} then 
$M(u)=\sum_w \mu(u,w) w$, but we will not need this fact to do our
computations.  Indeed, in Section~\ref{rfm} 
we will use the displayed definitions of $Z(u)$
and $M(u)$ above to reprove Theorem~\ref{muP*}.

The crucial observation underlying our method is that $Z(u)$ and
$M(u)$ can be expressed in terms of simpler series.  To define these
series, it will help to have a bit more notation.  If $S\sbe[\bmn]^*$
then we will also let $S$ stand for the generating function
$\sum_{w\in S} w$.  Context will make it clear which interpretation is
meant.  If $S$ is empty then the corresponding generating function is
the zero series.
If $f$ is a series without constant term then we let
$$
f^+=f+f^2+f^3+\cdots =f^*-\ep.
$$
Note that $f^+$ is rational if $f$ is.
Finally, a function $F:[\bmn]^*\ra \bbZ\sspn{[\bmn]}$ is called 
{\it multiplicative\/} if for any $u\in[\bmn]^*$ we have
$$
F(u)=F(u(1))F(u(2))\cdots F(u(l))
$$
where $l=\ell(u)$.

Now define two multiplicative functions from $[\bmn]^*$ to
$\bbZ\sspn{[\bmn]}$ by setting, for all $\bmk\in[\bmn]$,
$$
z(\bmk)=[\bmk,\bmn]\cdot[\bmk\bmmi\1]^*
$$
and
$$
m(\bmk)=\case{\1-\2^+(\ep-\1)}{if $k=1$,}
{\left(\bmk^+-(\bmk\bmpl\1)^+\right)(\ep-\1)}{if $k\ge2$.\rp{0}{15}}
$$
(Note that by convention, $[\bmk\bmmi\1]=\emp$ when $k=1$ 
and $\bmk\bmpl\1=\emp$ when $k=n$.) These are the building blocks for $Z(u)$ and
$M(u)$.
\ble
\label{ZM}
For any $u\in[\bmn]^*$ we have
$$
Z(u)=[\bmn]^* z(u)
$$
and
$$
M(u)=(\ep-\1) m(u)
$$
\ele
\begin{proof}
To prove the first equation, it suffices to show that 
the product on the right-hand side produces each $w\ge u$ according to
the rightmost embedding $\rho_u$ of $u$ into $w$.  So such $w$ will occur
exactly once since the rightmost embedding is unique.  Suppose $\bmk$
is the last element of $u$.  Then $z(\bmk)$ is the last factor of the
product.  The term $\bml$ chosen from $[\bmk,\bmn]$ corresponds to the
element of $w$ greater than $\bmk$ in the rightmost embedding, while
the product $[\bmk\bmmi\1]^*$ contains all possible subwords which
could appear after $\bml$ in $w$ while keeping $\bmk$ in its
rightmost position.  Similar considerations apply to the other factors
in $z(u)$.  Finally, the initial $[\bmn]^*$ accounts for everything to
the left of the element of $w$ corresponding to the first element of $u$.

The proof of the second equation is similar except that we must have a
unique term for every normal embedding $\eta_u$ into $w$ and each term
must have sign $(-1)^{d(\eta_u)}$.  Again, consider the last element
$\bmk$ of $u$.  If $\bmk=\1$ then by the first normality condition,
the corresponding element of $w$ must be $\bml=\1$ or $\bml=\2$.  If
$\bml=\1$ then the second normality condition ensures that there is
no element to the right of $\bml$ in $w$ and there is no contribution to the
defect in this case.  This corresponds to the initial $\1$ in the
expression for 
$m(\1)$.  If $\bml=\2$ then (by normality again) the subword of $w$ to
the right of $\bml$ must consist only of $\2$'s, possibly with a final
$\1$.  The factor $-\2^+$ accounts for the string of $\2$'s with the
appropriate sign and the final factor of $\ep-\1$ takes care of the
possibilities at the right end of $w$.  The arguments for $\bmk\ge\2$ 
and for the initial factor of $\ep-\1$ in $M(u)$ are along the same lines
and so omitted. 
\end{proof}

Note that directly from their definition, $z(\bmk)$ and $m(\bmk)$ are
rational series.  So, by the previous lemma, we have the following result.
\bth
For any $u\in[\bmn]^*$, $Z(u)$ and $M(u)$ are rational series.\qed
\eth

We will now prove analogous results for the generating functions of
$\ze$ and $\mu$ using the alphabet of ordered pairs
$[\bmn]\times[\bmn]=[\bmn]^2$.   We could do so by modifying the
arguments which led to the previous theorem.  But for variety's sake,
we will use finite-state automata.  We write the elements of 
$\bbZ\sspn{[\bmn]^2}$ as
$$
f=\sum_{u,w} c_{u,w}\ u\otimes w.
$$

Given an alphabet $A$, a {\it finite-state automaton\/} is a digraph $D$
with the following properties.  The vertex set $V$ and directed edge
(arc) set $E$ are both finite with loops and multiarcs permitted.  There
is a distinguished initial vertex and a distinguished final vertex
denoted $\al$ and $\om$, respectively.  Each $e\in E$ is assigned a
monomial label $f(e)\in\bbZ\sspn{A}$.

Now given a finite walk $W$ with arcs $e_1,\ldots,e_l$, we assign it the
monomial 
$$
f(W)=\prod_{i=1}^l f(e_i).
$$
The formal power series {\it accepted\/} by $D$ is
$$
f(D)=\sum_W f(W),
$$
where the sum is over all finite walks from $\al$ to $\om$.  Note that
if $e_1,\ldots,e_j$ are all arcs from a vertex $\be$ to a vertex $\ga$,
then replacing these arcs by a single arc $e=\widevec{\be\ga}$ and
setting
$$
f(e)=\sum_{i=1}^j f(e_i)
$$
does not change the series accepted by $D$.  So we will do this when
constructing automata without further comment.  We will also use algebraic
operations to simplify the sum for $f(e)$ if possible.  

The crucial
fact which we will need is the well-known result that a series  is
rational if and only if it is accepted by some finite-state automaton $D$,
see e.g. \cite{br:rsl}.

\bfi
$$
\begin{psmatrix}[rowsep=3cm,colsep=4.5cm]
                        &\pscb{\rnode{al}{\al}}  &\\
\pscb{\rnode{1}{\be_1}} &\pscb{\rnode{2}{\be_2}} &\pscb{\rnode{3}{\be_3}}
\end{psmatrix}
\psset{nodesep=5pt,arrows=->}
\everypsbox{\scriptstyle}
\ncline{al}{1}\Bput{\1\otimes(\1+\2+\3)}
\ncline[nodesep=8pt]{al}{2}\Aput{\2\otimes(\2+\3)}
\ncline{al}{3}\Aput{\3\otimes\3}
\nccircle[angleA=0,nodesepA=8pt]{al}{.5}\Bput{\ep\otimes(\1+\2+\3)}
\nccircle[angleA=90,nodesepA=8pt]{1}{.5}\Bput{\1\otimes(\1+\2+\3)}
\nccircle[angleA=180,nodesepA=8pt]{2}{.5}\Bput{\ep\otimes\1+\2\otimes(\2+\3)}
\nccircle[angleA=270,nodesepA=8pt]{3}{.5}\Bput{\ep\otimes(\1+\2)+\3\otimes\3}
\ncarc[nodesep=8pt]{1}{2}\Aput{\2\otimes(\2+\3)}
\ncarc[nodesep=8pt]{2}{1}\Aput{\1\otimes(\1+\2+\3)}
\ncarc[nodesep=8pt]{2}{3}\Aput{\3\otimes\3}
\ncarc[nodesep=8pt]{3}{2}\Aput{\2\otimes(\2+\3)}
\ncangle[angleA=-90,angleB=-90,arm=2,linearc=.15,nodesep=8pt]{1}{3}\Bput{\3\otimes\3}
\ncangle[angleA=90,angleB=90,arm=5.5,linearc=.15,nodesep=8pt]{3}{1}\Bput{\1\otimes(\1+\2+\3)}
$$
\vspace{2cm}
\capt{The automaton for $Z_\otimes$ when $n=3$}
\label{aZ}
\efi

\bth
In $\bbZ\sspn{[\bmn]^2}$ the series
$$
Z_\otimes=\sum_{u,w} \zeta(u,w) u\otimes w
$$
and
$$
M_\otimes=\sum_{u,w} \mu(u,w) u\otimes w
$$
are rational.
\eth
\begin{proof}
For both series, we will build finite-state automata accepting them.

The automaton $D$ for $Z_\otimes$ has vertices
$\{\al,\om,\be_1,\ldots,\be_n\}$.  
A picture of the digraph when $n=3$ is given in Figure~\ref{aZ}.
The vertex $\om$ is not shown since it simply has an incoming  arc,
labeled $\ep\otimes\ep$, from every other vertex.
To describe the arc set, we will consider each of the other 
vertices in turn and describe all its incoming arcs.

If the vertex is $\al$, then the only incoming arc is a loop labeled
$\ep\otimes[\bmn]$.  If it is $\om$, then we have already described the
arcs into it.  If the vertex is $\be_k$ for some $k$ then
there is an incoming arc from every vertex except $\om$,
as well as a loop, which are labeled 
$$
f(\widevec{\be\be_k})=
\case{\ep\otimes[\bmk\bmmi\1]+\bmk\otimes[\bmk,\bmn]}{if $\be=\be_k$,}
{\bmk\otimes[\bmk,\bmn]}{else.\rp{0}{15}}
$$

To show $D$ accepts $Z_\otimes$, we need to prove that for every pair 
$u\otimes w$ with $u\le w$ there is a unique way to obtain 
$u\otimes w$ as a monomial along some walk from $\al$ to $\om$, and
that these are the only monomials in $f(D)$.  We will indicate how one
can find the walk $W$ given $u\otimes w$, since then the reader should be
able to fill in the details of the rest of the proof.  In fact, we
will show that $W$ constructs $w$ and $u$ in its rightmost embedding $\rho_u$
into $w$ in the following sense.  If $e_i$ is the $i$th arc of $W$
then $f(e_i)$ contains the term $a\otimes b$ where $b=w(i)$ and
$a=\rho_u(i)$ or $\ep$ depending on whether $\rho_u(i)\in[\bmn]$ or
$\rho_u(i)=\0$, respectively.

To begin, $W$ loops $i-1$ times at $\al$,  where $i$ is the smallest index
with $\rho_u(i)\neq\0$. (If $u=\ep$ then let $i=\ell(w)+1$.)  The walk $W$
finishes at $\om$ if $i=\ell(w)+1$, while if $i\le\ell(w)$ it goes to 
$\be_k$ where $\rho_u(i)=\bmk$.  
Now $W$ loops at $\be_k$ through arc
$e_{j-1}$, where $j>i$ is the next index with $\rho_u(j)\neq0$.  
The $\ep\otimes[\bmk\bmmi\1]$ summand on the arc contains the
necessary monomial.  Then $e_j$ goes from $\be_k$ to $\be_l$ where
$\rho_u(j)=\bml$.  Note that we could have $k=l$ so that this would also
be a loop, in which case the $[\bmk,\bmn]$ summand contains the
desired monomial.  One continues in this manner until $W$ has gone
through $\ell(w)$ arcs, after which it takes the arc to $\om$.

The automaton for $M_\otimes$ has the same vertex set as the one for
$Z_\otimes$.  See Figure~\ref{aM} for the picture when $n=3$.
Again, $\om$ only has incoming arcs from the other vertices and
they are all labeled $\ep\otimes\ep$, so it is not shown.  Since the
construction of this automaton and the proof that it does accept $M_\otimes$ is
parallel to what we did for $Z_\otimes$, we will content ourselves with a
description of its arc set. Note that the interpretation of $\mu(u,w)$
built into the automaton relies on Theorem \ref{muP*}.

For $\al$ there are no incoming arcs and we have already
described what happens for $\om$.  If  the vertex is $\be_1$, then
there are incoming arcs from every vertex except $\om$ and they are labeled
$$
f(\widevec{\be\be_1})=
\case{\1\otimes\1}{if $\be=\be_1$,}{(\1-\ep)\otimes\1}{else.\rp{0}{15}}
$$
If the vertex is $\be_k$ for $k\ge2$ then we have the same set of
incoming arcs with labels
$$
f(\widevec{\be\be_k})=
\case{(\bmk-(\bmk\bmmi\1)+\ep)\otimes\bmk}
{if $\be=\be_k$,}
{\bmk\otimes\bmk}
{else.\rp{0}{15}}
$$
This completes the description of the automaton for $M_\otimes$.
\end{proof}

\bfi
$$
\begin{psmatrix}[rowsep=3cm,colsep=4.5cm]
                        &\pscb{\rnode{al}{\al}}  &\\
\pscb{\rnode{1}{\be_1}} &\pscb{\rnode{2}{\be_2}} &\pscb{\rnode{3}{\be_3}}
\end{psmatrix}
\psset{nodesep=5pt,arrows=->}
\everypsbox{\scriptstyle}
\ncline{al}{1}\Bput{(\1-\ep)\otimes\1}
\ncline[nodesep=8pt]{al}{2}\Aput{\2\otimes\2}
\ncline{al}{3}\Aput{\3\otimes\3}
\nccircle[angleA=90,nodesepA=8pt]{1}{.5}\Bput{\1\otimes\1}
\nccircle[angleA=180,nodesepA=8pt]{2}{.5}\Bput{(\2-\1+\ep)\otimes\2}
\nccircle[angleA=270,nodesepA=8pt]{3}{.5}\Bput{(\3-\2+\ep)\otimes\3}
\ncarc[nodesep=8pt]{1}{2}\Aput{\2\otimes\2}
\ncarc[nodesep=8pt]{2}{1}\Aput{(\1-\ep)\otimes\1}
\ncarc[nodesep=8pt]{2}{3}\Aput{\3\otimes\3}
\ncarc[nodesep=8pt]{3}{2}\Aput{\2\otimes\2}
\ncangle[angleA=-90,angleB=-90,arm=2,linearc=.15,nodesep=8pt]{1}{3}\Bput{\3\otimes\3}
\ncangle[angleA=90,angleB=90,arm=4.5,linearc=.15,nodesep=8pt]{3}{1}\Bput{(\1-\ep)\otimes\1}
$$
\vspace{2cm}
\capt{The automaton for $M_\otimes$ when $n=3$}
\label{aM}
\efi

\section{Generating functions in commuting variables}
\label{gfc}

By specialization of variables, we can get generating functions for
$\ze$ and $\mu$ it terms of the length function $\ell(w)$ or in terms of the sum
of the parts, or {\it norm}, of the composition, which will be denoted $|w|$.  We will
also need to keep track of the {\it type of $w$},
$t(w)=(l_1,l_2,\ldots,l_n)$, where $l_k$ is the number of $\bmk$'s in
$w$.  So $\sum_k l_k =\ell(w)$ and $\sum_k l_k k=|w|$.

Suppose $x$ is a variable and we substitute $x^k$ for $\bmk$ in
$Z(u)$.  Then the generating function becomes 
$$
Z(u;x)=\sum_{w\ge u}x^{|w|}.
$$
Doing the same thing with $z(\bmk)$ and summing the resulting
geometric series gives
$$
z(\bmk;x)=\frac{x^k+x^{k+1}+\cdots+x^n}{1-(x+x^2+\cdots+x^{k-1})}
=\frac{x^k-x^{n+1}}{1-2x+x^k}.
$$
If $t(u)=(l_1,\ldots,l_n)$, then appealing to Lemma~\ref{ZM} yields a
norm generating function in $[\bmn]^*$ of
$$
Z(u;x)=\frac{1-x}{1-2x+x^{n+1}}
    \prod_{k=1}^n \left(\frac{x^k-x^{n+1}}{1-2x+x^k}\right)^{l_k}.
$$
Note that this generating function depends only on the type of $u$ and
not on $u$ itself.
Note also that one can take $n\ra\infty$ in this series (reflecting the
fact that there are only finitely many compositions with given norm)
to obtain the norm generating function in $\bbP^*$
$$
Z_\bbP(u;x)=\frac{1-x}{1-2x}
    \prod_{k\ge1} \left(\frac{x^k}{1-2x+x^k}\right)^{l_k}.
$$
When $u=\ep$, this shows that the rank generating function for 
$\bbP^*$ (which is graded by norm) is $(1-x)/(1-2x)$. This can also be
seen from the fact that there are $2^{N-1}$ compositions of $N$ for $N\ge1$.
This same procedure can be applied to the generating function $M(u)$.

If one wants the generating function by length, then one substitutes
the same variable, say $t$, for each $\bmk$.  
Under this substitution
$m(\bmk;t)=0$ for $1\le k\le n$ and so $M(u;t)=0$ unless $u=\ep$.
Also, in this case one needs
to remain in $[\bmn]^*$ since there are infinitely many compositions
in $\bbP^*$ of a given nonzero length.  The details of these
computations are routine, so we will merely state
the results.
\bth
Let $t(u)=(l_1,\ldots,l_n)$ where $u\in[\bmn]^*$.  Then we have the
norm generating functions
$$
Z(u;x)=\frac{1-x}{1-2x+x^{n+1}}
    \prod_{k=1}^n \left(\frac{x^k-x^{n+1}}{1-2x+x^k}\right)^{l_k}
$$
and
$$
M(u;x)=\frac{x^{|u|}(1-x)^{2\ell(u)+1}}{(1-x)^{l_1+l_n}}
    \prod_{k=2}^n \frac{1}{(1-x^k)^{l_{k-1}+l_k}}.
$$
We also have the length generating functions
$$
Z(u;t)=\frac{1}{1-nt}
    \prod_{k=1}^n \left(\frac{(n-k+1)t}{1-(k-1)t}\right)^{l_k}
$$
and
$$
M(u;t)=\case{1-t}{if $u=\ep$,}{0}{else.}
$$

In $\bbP^*$ we have norm generating functions
$$
Z_\bbP(u;x)=\frac{1-x}{1-2x}
    \prod_{k\ge1} \left(\frac{x^k}{1-2x+x^k}\right)^{l_k}
$$
and\\
\eqed{\dil
M_\bbP(u;x)=\frac{x^{|u|}(1-x)^{2\ell(u)+1}}{(1-x)^{l_1}}
    \prod_{k\ge2} \frac{1}{(1-x^k)^{l_{k-1}+l_k}}.
}
\eth

We would now like to calculate the generating function for $\ze^m$.
This is of interest because $\ze^m(u,w)$ counts the number of
multichains of length $m$ from $u$ to $w$.  (As mentioned in the
introduction, the original motivation of Bergeron et.\ al.\ in
studying $\bbP^*$ was to count saturated chains in $[\ep,w]$.)  To do
this, we will have to exploit a connection between the
incidence algebra $I([\bmn]^*)$ and the algebra
$\End\bbZ\sspn{[\bmn]}$ of continuous linear endomorphisms of
$\bbZ\sspn{[\bmn]}$ (for the meaning of ``continuity'' here, see
e.g. \cite[p. 55]{br:rsl}).  This relationship will also be important in the
next section where we will reprove the formula for $\mu$.

Note that~\ree{Zu} already defines a map $Z:[\bmn]^*\ra \bbZ\sspn{[\bmn]}$.  
We can extend this to an element of $\End\bbZ\sspn{[\bmn]}$ as
follows.  Take any $\phi\in I([\bmn]^*)$ and define a corresponding map
$F_\phi:[\bmn]^*\ra \bbZ\sspn{[\bmn]}$ by
$$
F_\phi(u)=\sum_w \phi(u,w)w
$$
where the sum is over all $w\in[\bmn]^*$, or equivalently over all
$w\ge u$ since $\phi(u,w)=0$ otherwise.  By continuity and linearity,
we can extend $F_\phi$ to a function in $\End\bbZ\sspn{[\bmn]}$ by letting
$$
F_\phi\left(\sum_u c_u u\right)=\sum_u c_u F_\phi(u).
$$
Note that the right-hand side converges since any $v\in[\bmn]^*$ 
occurs with nonzero coefficient in only finitely many of the summands
$F_\phi(u)$.  Lifting elements of $I([\bmn]^*)$ to $\End\bbZ\sspn{[\bmn]}$
in this way is well behaved.

\bth
\label{iso}
The map $\phi\mapsto F_\phi$ is an algebra anti-isomorphism of 
$I([\bmn]^*)$ with a subalgebra of $\End\bbZ\sspn{[\bmn]}$.
\eth
\begin{proof}
Checking the various needed properties of the map are easy, so we will
just indicate why multiplication is antipreserved to illustrate.  Recall
that the product of $\phi,\psi\in I([\bmn]^*)$ is their convolution
$\phi*\psi$ while the product in $\End\bbZ\sspn{[\bmn]}$ is
composition of functions.  To show that the two multiplications
correspond, it suffices to check that they do so on elements
$u\in[\bmn]^*$.   So we compute
\bea
F_\psi\circ F_\phi(u)&=&F_\psi\left(\sum_v \phi(u,v) v\right)\\
      &=&\sum_{v,w} \phi(u,v)\psi(v,w)w\\
      &=&\sum_w \phi*\psi(u,w)w\\
      &=&F_{\phi*\psi}(u)
\eea
as desired.
\end{proof}

Now we can factor the generating function for $\ze^m$ as follows.  Let 
$\bbZ[[X]]$ be the formal power series ring over the integers in the
set $X=\{x_1,x_2,\ldots,x_n\}$ of commuting variables.  Consider the
projection map $\rho:\bbZ\sspn{[\bmn]}\ra\bbZ[[X]]$ which sends $\bmk$
to $x_k$.  Then we have
$$
\rho\circ z(\bmk)=\frac{x_k+\cdots+x_n}{1-x_1-\cdots-x_{k-1}}.
$$
Define a multiplicative function $f:\bbZ[[X]]\ra\bbZ[[X]]$ by
\beq
\label{f}
f(x_k)=\frac{x_k+\cdots+x_n}{1-x_1-\cdots-x_{k-1}}.
\eeq
Clearly $f$ is constructed so that
$$
\rho\circ z = f\circ \rho.
$$

We now apply the same idea to the function $Z$.  
If $u\in[\bmn]^*$ then we let $X^u=\prod_k x_k^{l_k}$, where
$t(u)=(l_1,\ldots,l_n)$.  So $\rho(u)=X^u$. 
Define a continuous, linear map $F:\bbZ[[X]]\ra\bbZ[[X]]$ by
$$
F(X^u)=\frac{1}{1-x_1-\cdots-x_n}f(X^u).
$$
It follows that
$$
\rho\circ Z = F\circ \rho.
$$

From Theorem~\ref{iso} we have that 
$$
\sum_w \ze^m(u,w) w =Z^m(u).
$$
So letting $t(u)=(l_1,\ldots,l_n)$ and applying $\rho$ to both sides,
we see that the generating function for $\ze^m$ in  $\bbZ[[X]]$ is
\bea
\sum_w \ze^m(u,w) X^w 
  &=&\rho\circ Z^m(u)\\
  &=&F^m\circ \rho(u)\\
  &=&F^m(X^u)\\
  &=&\prod_{i=0}^{m-1}\frac{1}{1-f^i(x_1)-\cdots-f^i(x_n)}
     \prod_{k=1}^n \left(f^m(x_k)\right)^{l_k}
\eea
where the last equality follows from an easy induction on $m$.

Thus to find $\ze^m$ for all $m$, it suffices to find $f^m(x_k)$ for all
$m$ and $k$.  Since this turns out to be surprisingly hard to do, we
will just consider what happens when $n=2$.  This case is of independent
interest because then the poset has rank numbers given by the
Fibonacci sequence.  However, this is different from the Fibonacci
posets defined by Stanley~\cite{sta:fl,sta:dp}.

For simplicity when $n=2$, let $x=x_1$ and $y=x_2$.  In this
case~\ree{f} becomes
$$
f(x)=x+y\qmq{and} f(y)=\frac{y}{1-x}.
$$
To simplify notation again, let 
$$
a_m=f^m(x)\qmq{and} b_m=f^m(y).
$$
Now we have, for $m\ge1$.
$$
f^m(x)=f^{m-1}(f(x))=f^{m-1}(x+y)=f^{m-1}(x)+f^{m-1}(y)
$$
or
\beq
\label{am}
a_m=a_{m-1}+b_{m-1}.
\eeq
Similarly, one obtains
\beq
\label{bm}
b_m=\frac{b_{m-1}}{1-a_{m-1}}
\eeq
for $m\ge1$, and it is easy to see that
\beq
\label{a0}
a_0=x\qmq{and} b_0=y.
\eeq
Hence we have to solve two recurrence relations in two unknowns.

Let us first make the norm substitution $y=x^2$.  In this case we will
denote $a_m$ and $b_m$ by $a_m(x)$ and $b_m(x)$.  To state our result,
we will need the round-down function $\fl{\cdot}$ and round-up
function $\ce{\cdot}$.  We will also use the conventions that the
binomial coefficient ${n\choose k}$ equals $0$ for $k<0$ or $k>n$ and
equals $1$ for $k=0$ and any $n$.
\bth
\label{zemx}
Suppose $u\in[\2]^*$ has type $t(u)=(l_1,l_2)$.  Then 
$$
\sum_{w} \ze^m(u,w) x^{|w|}= a_m(x)^{l_1} b_m(x)^{l_2}
\prod_{i=0}^{m-1}\frac{1}{1-a_i(x)-b_i(x)}.
$$
Furthermore, for all $m\ge0$ we have
\beq
\label{amx}
a_m(x)=\frac{x\ab_m(x)}{d_m(x)}\qmq{and} b_m(x)=\frac{x^2}{d_m(x)d_{m+1}(x)}
\eeq
where
\beq
\label{abmx}
\ab_m(x)=\sum_i (-1)^{\flf{i}{2}}{\flf{m+i}{2}\choose i} x^i
\qmq{and}
d_m(x)=\sum_i (-1)^{\cef{i}{2}}{\flf{m+i-1}{2}\choose i} x^i.
\eeq
\eth
\bprf
It suffices to show that the equations for $a_m(x)$ and $b_m(x)$ given
in the statement of the 
theorem satisfy~\ree{am}, \ree{bm}, and~\ree{a0}.  Checking the
boundary conditions is easy.

To prove that~\ree{am} holds, substitute~\ree{amx} into the recursion,
multiply by $d_m(x)d_{m+1}(x)/x$, substitute~\ree{abmx}, and take the
coefficient of $x^k$ on both sides.  Thus we need to prove
\beq
\label{sum}
\barr{l}
\dil\sum_i (-1)^{\flf{i}{2}+\cef{k-i}{2}}
{\flf{m+i}{2}\choose i}{\flf{m+k-i-2}{2}\choose k-i}\\[20pt]
\qquad =
\dil\sum_i (-1)^{\flf{i}{2}+\cef{k-i}{2}}
{\flf{m+i-1}{2}\choose i}{\flf{m+k-i-1}{2}\choose k-i}
\earr
\eeq
for $k\neq1$.  (When $k=1$ we need to add a 1 onto the right-hand side
corresponding to the $x$ obtained from $b_m(x)$ after doing the
multiplication.  But this identity is easy to verify.)
The proof now breaks down into four cases depending on the parities of
$m$ and $k$.  We will only discuss what happens when $m$ is even and $k$
odd, as the other demonstrations are similar.

So suppose $m=2l$ and $k=2j+1$ for integers $l,j$.  Then the terms
in~\ree{sum} corresponding to even $i$ cancel.  Rewriting the odd $i$
terms using rising factorials yields, after some cancellation, the
equivalent hypergeometric series identity
$$
\barr{l}
(l-1)_{j+1}(2-l)_j 
\hspace{5pt}
\rule{0in}{0in}_4F_3
\left[\barr{ccccc}
  l+1,   &-l,    &-j,    &-j-1/2;    &1\\
         &l+j-1, &1-l-j  &1/2        &
\earr
\right]
\\[20pt]
\qquad =
(l)_{j+1}(1-l)_j
\hspace{5pt}
\rule{0in}{0in}_4F_3
\left[\barr{ccccc}
  l,     &1-l,    &-j,    &-j-1/2;    &1\\
         &l-j,    &-l-j   &1/2        &
\earr
\right].
\earr    
$$
Using the implementation of Zeilberger's
algorithm~\cite{zei:fap,zei:hsa} due to Paule and
Schorn~\cite{ps:mvz}, one can verify that both sides of this equation
satisfy the same three-term recurrence relation in $l$.  Also, $j\neq0$
since $k\neq 1$.  For positive $j$ both sides of the equation are
clearly zero for $l=0,1$.  So since both sides also satisfy the same
boundary conditions, they must be equal.  Also, Dennis Stanton has 
pointed out that one can give a more traditional proof of this
identity (and, in fact, prove a generalization of it) using Tchebyshev
polynomials and trigonometric identities. 

Verifying~\ree{bm} turns out to be much simpler.
Substituting~\ree{amx}, clearing denominators, and dividing by $x^2$,
leads to the equivalent identity 
$$
x\ab_{m-1}(x)+d_{m+1}(x)-d_{m-1}(x)=0.
$$
This follows easily from~\ree{abmx} and the binomial recursion.
\eprf

To get the corresponding length generating functions, we need only
change the boundary conditions to $x=y=t$.  In this case we write
$a_m(t)$ and $b_m(t)$ for $a_m$ and $b_m$.  Since the computations are
similar, we will simply state the result.
\bth
\label{zemt}
Suppose $u\in[\2]^*$ has type $t(u)=(l_1,l_2)$.  Then 
$$
\sum_{w} \ze^m(u,w) t^{\ell(w)}= a_m(t)^{l_1} b_m(t)^{l_2}
\prod_{i=0}^{m-1}\frac{1}{1-a_i(t)-b_i(t)}.
$$
Furthermore, for all $m\ge0$ we have
$$
a_m(t)=\frac{t\ab_m(t)}{d_m(t)}\qmq{and} b_m(t)=\frac{t}{d_m(t)d_{m+1}(t)}
$$
where $\ab_m(t)=\sum_i (-1)^i \al_{m,i} t^i$
and $d_m(t)=\sum_i (-1)^i \de_{m,i} t^i$ with the coefficients 
$\al_{m,i}$ and $\de_{m,i}$ being given by
$$
\al_{m,i}=
\case{\dil \frac{(m+1)2^i}{2i+1}{\frac{m+2i}{2}\choose\frac{m-2i}{2}}}
{if $m$ is even,}
{\dil\rule{0pt}{30pt} 2^{i+1}{\frac{m+2i+1}{2}\choose\frac{m-2i-1}{2}}}
{if $m$ is odd,}
$$
and\\
\eqed{
\de_{m,i}=
\case{\dil \frac{m2^i}{m+2i}{\frac{m+2i}{2}\choose\frac{m-2i}{2}}}
{if $m$ is even,}
{\dil\rule{0pt}{30pt} 2^i{\frac{m+2i-1}{2}\choose\frac{m-2i-1}{2}}}
{if $m$ is odd.}
}
\eth

\section{Reproving the formula for $\mu$ in $\bbP^*$}
\label{rfm}
We will now reprove the formula for $\mu$ in Theorem~\ref{muP*}.  Our
principal tools will be the descriptions of $Z$ and $M$ in
Lemma~\ref{ZM} and the  anti-isomorphism in Theorem~\ref{iso}.
Although we only stated the latter result for $[\bmn]^*$, it clearly 
holds also for $\bbP^*$.  The Lemma must be modified slightly by letting
$n$ tend to $\infty$.  So the formulas for $z$ and $Z$ become
$$
z(\bmk)=[\bmk,\infty)[\bmk\bmmi\1]^*
$$
and
$$
Z(u)=\bbP^*z(u).
$$

{\it Proof (of Theorem~\ref{muP*}).}
We wish to show that $\zeta*\mu$ is the identity element of the
incidence algebra.  So by Theorem~\ref{iso}, it suffices to show that
$M\circ Z$ is the identity endomorphism.  For any $u\in\bbP^*$ we
have, using the multiplicativity of $m$,
$$
M\circ Z(u)= M(\bbP^*z(u))=(\ep-\1) m(\bbP)^* m(z(u)).
$$
So it will be enough to show
$$
(\ep-\1)m(\bbP)^*=\ep\qmq{and} m(z(\bmk))=\bmk
$$
for all $\bmk\in\bbP$.

For the first equation, note that
\bea
m(\bbP)&=&m(\1)+m(\2)+m(\3)+\cdots\\[3pt]
  &=&(\1-\2^+(\ep-\1))+(\2^+(\ep-\1)-\3^+(\ep-\1))+(\3^+(\ep-\1)-\4^+(\ep-\1))+\cdots\\[3pt]
  &=&\1.
\eea
Now~\ree{*} gives
$$
(\ep-\1)m(\bbP)^*=(\ep-\1)\1^*=(\ep-\1)(\ep-\1)^{-1}=\ep.
$$

For the second equation, we note that the case $\bmk=\1$ has already
been done in the previous paragraph, since
$$
m(z(\1))=m(\bbP)=\1.
$$
Note that for $k\ge2$ the same telescoping phenomenon gives 
$$
m([\bmk,\infty))=\bmk^+(\ep-\1)
\qmq{and}
m([\bmk\bmmi\1])=\1-\bmk^+(\ep-\1).
$$
Combining this with~\ree{*}, we obtain
\bea
m(z(\bmk))&=&m([\bmk,\infty))m([\bmk\bmmi\1])^*\\[3pt]  
  &=&\bmk^+(\ep-\1)\left(\1-\bmk^+(\ep-\1)\right)^*\\[3pt]
  &=&\bmk^+(\ep-\1)\left(\ep-\1+\bmk^+(\ep-\1)\right)^{-1}\\[3pt]
  &=&\bmk^+(\ep-\1)\left((\ep+\bmk^+)(\ep-\1)\right)^{-1}\\[3pt]
  &=&\bmk^+(\ep-\1)(\ep-\1)^{-1}(\ep+\bmk^+)^{-1}\\[3pt]
  &=&\bmk\bmk^*\left(\bmk^*\right)^{-1}\\[3pt]
  &=&\bmk.
\eea
This finishes the proof of Theorem~\ref{muP*}.
\qed

\section{Generalized subword order}
\label{gso}

We now present a rubric due to Sagan and Vatter~\cite{sv:mfc} under
which the theorems about rationality of the M\"obius and zeta
functions for $A^*$ and $\bbP^*$ both become special cases.  Let $P$
be any poset. Turn $P^*$ into a poset by letting $u\le_{p^*} w$ if
there is a subword $w(i_1)\ldots w(i_l)$ of $w$ having length $l=\ell(u)$
such that
$$
\mbox{$u(i)\le_P w(i_l)$ for $1\le i\le l$.}
$$
We call this the {\it generalized subword order\/} on $P^*$.  Note
that we recover $A^*$ or $\bbP^*$ if we take $P$
to be an antichain or a well-ordered countably infinite chain, respectively.
Note also that we will leave off the subscripts on inequalities if it
is clear from context which poset is meant.

Many of our results about $\ze$  for $\bbP^*$ from Sections~\ref{r} and~\ref{gfc},
as well as the corresponding ones for $A^*$ 
of Bj\"orner and Reutenauer~\cite{br:rsl},
generalize easily to $P^*$.  Given an element $a\in P$ we consider the
{\it upper order ideal generated by $a$} and its set-theoretic complement
$$
\mbox{$I_a=\{c\in P\ |\ c\ge_P a\}$ and $J_a=P-I_a$,}
$$
respectively.  We define $Z(u)$ in $P^*$ by~\ree{Zu} as before and
also define a multiplicative map from $P^*$ to $\bbZ\sspn{P}$ by
$$
z(a)=I_a J_a^*.
$$
The proofs we have already seen contain all the ideas needed to
demonstrate the next result, so we suppress the details.  We will also
use the same notation as in the earlier results, as we did with $Z(u)$.
\bth
Let $P$ be any poset.  Then for any $u\in P^*$ we have
$$
Z(u)=P^*z(u)
$$
and so $Z(u)$ is rational.  Similarly, in $\bbZ\sspn{P^2}$ the series
$$
Z_\otimes =\sum_{u,w}\zeta(u,w)u\otimes w
$$
is rational.  Finally, if $P$ is finite and
$u$ has $l_a$ occurrences of $a$ for each $a\in P$, then we have the
length generating function\\
\eqed{
Z(u;t)=\frac{1}{1-|P|t}\prod_{a\in P}
\left(\frac{|I_a|t}{1-|J_a|t}\right).
}
\eth

Generalizing our results about $\mu$ is more delicate.  Indeed, there
is no known formula for the M\"obius function in $P^*$ for arbitrary
$P$.  However, there is a class of posets for which $\mu$ has been
found.  To characterize the M\"obius function in these posets, we need
the appropriate definition of a normal embedding.
Suppose $\zh$ is a new element not in $P$ and form a poset
$\Ph$ on $P\cup\zh$ by adding the relations $\zh<_{\Ph} a$ for all 
$a\in P$.  One defines support and expansion exactly as before, just
replacing $0$ with $\zh$.  Then for $u,w\in P^*$, an {\it embedding\/}
of $u$ into $w$ 
is an expansion $\eta_u\in\Ph^*$ of length $\ell(w)$ such that 
$$
\mbox{$\eta_u(i)\le_{\Ph} w(i)$ for $1\le i\le \ell(w)$.}
$$
Clearly, $u\le_{P^*} w$ if and only if there is an embedding of
$u$ into $w$.

To define normality, call $P$ a {\it rooted tree\/} if its Hasse
diagram is a tree having a unique minimal element.  More generally,
call $P$ a {\it rooted forest\/} if the connected components
of its Hasse diagram are rooted trees.  Note that in this case $\Ph$
is a rooted tree.  So given $a\in P$ we can define $a^-$ to be the
element adjacent to $a$ on the unique path from $a$ to $\zh$ in $\Ph$.
If $P$ is a rooted forest, define an embedding $\eta_u$ of $u$ into $w$ to be
{\it  normal\/} if it satisfies the following pair of conditions.
\begin{enumerate}
\item  For $1\le i\le \ell(w)$ we have $\eta_u(i)=w(i)$,
$w(i)^-$, or $\zh$.\vs{3pt}
\item  For all $a\in P$ and every run $[r,t]$ of $a$'s
in $w$, we have\vs{3pt}
\begin{enumerate}
\item $(r,t]\sbe \Supp \eta_u$ if  $a$ is minimal in $P$,\vs{3pt}
\item $r\in\Supp \eta_u$ otherwise.
\end{enumerate}
\end{enumerate}

In this situation, the definition of the {\it defect\/} of a normal
embedding $\eta_u$ into $w$ should come as no surprise:
$$
d(\eta_u)=\#\{i\ |\ \eta_u(i)=w(i)^-\}.
$$
The following theorem generalizes both Theorem~\ref{muA*} and
Theorem~\ref{muP*}.
\bth[Sagan and Vatter~\cite{sv:mfc}]
\label{mufor}
Let $P$ be a rooted forest.  Then the M\"obius function of $P^\ast$ is
given by
$$
\mu(u,w)=\sum_{\eta_u}(-1)^{\eta_u},
$$
where the sum is over all normal embeddings $\eta_u$ of $u$ into
$w$.\qed
\eth

With this result in hand, generalizing the results for $\mu$ follows
the same lines as for $\ze$.  If $P$ is any poset then let $O_P$ be
the set of minimal elements of $P$.  (So if $P=\bbP$ then 
$O_P=\{1\}$.)  Also, if $a\in P$ then the set of elements {\it
covering\/} $a$ is
$$
C_a =\{c\in P\ |\ \mbox{$c>a$ and there is no $b$ with $c>b>a$}\}.
$$
Now let $P$ be a rooted forest and define $M(u)$ for
$u\in P^*$ by equation~\ree{Mu}.  The corresponding multiplicative
function is
$$
m(a)=
\case{\dil a-\left(\sum_{c\in C_a} c^+\right)(\ep-O_P)}
{if $a\in O_P$,}
{\dil\rule{0pt}{30pt} \left(a^+-\sum_{c\in C_a} c^+\right)(\ep-O_P)}
{else.}
$$

Again, there is nothing really new in considering an arbitrary rooted
forest instead of $\bbP$, so we will merely state the results.
\bth
Let $P$ be a rooted forest.  Then for any $u\in P^*$ we have
$$
M(u)=(\ep-O_P)z(u)
$$
and so $M(u)$ is rational.  Similarly, in $\bbZ\sspn{P^2}$ the series
$$
M_\otimes =\sum_{u,w}\mu(u,w)u\otimes w
$$
is rational.  Finally, if $P$ is finite
then we have the length generating function
$$
M(u;t)=\frac{t^{|P|}(1-|O_P|t)^{|P-O_P|+1}}{(1-t)^{|P|}}
\prod_{a\in O_P}\left(1-t-|C_a|(1-|O_P|t)\right)
\prod_{b\not\in O_P}(1-|C_b|).
$$
In particular, if $u$ contains any element which is covered by exactly one
other element then $M(u;t)=0$.\qed
\eth

As a final remark, one can give a proof of Theorem~\ref{mufor} in the
same way as was done for Theorem~\ref{muP*} in the previous section.

\section{Comments and open problems}
\label{cop}

We end with some comments and open problems.

\subsection{Generating functions for $\ze^m$}

It would be interesting to compute the generating function for $\ze^m$
in $[\bmn]^*$ for arbitrary $n$.  It appears that one
can say something, at least for $n=3$.  Let $a_m(t),b_m(t),c_m(t)$
stand for $f^m(x_1),f^m(x_2),f^m(x_3)$, respectively, when using the
length generating function.  Then numerical evidence suggests that
there is a polynomial $d_m(t)$ such that the denominators of our three
rational functions factor as $d_1d_2\cdots d_{2m-2}$,
$d_{2m-3}d_{2m-2}d_{2m-1}$, and $d_{2m-2}d_{2m}$, respectively.  Note
that the behaviour of the denominator of $a_m(t)$ behaves differently
from the $n=2$ case in that the number of factors increases with $m$.

It would also be interesting to find ``classical'' proofs of the
hypergeometric identities used in the demonstrations of
Theorems~\ref{zemx} and~\ref{zemt}.  The series involved are neither
0-balanced nor well-poised so we were unable to come up with appropriate
theorems in the literature which applied to them.  Andrew Sills has
noted that they are 1-balanced, which may be of help.

\begin{figure}
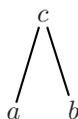

\begin{center}
\begin{footnotesize}
$$
\psset{nodesep=2pt,colsep=18pt,rowsep=25pt}
\begin{psmatrix}
[name=c] c\psspan{2}\\[0pt]
[name=a] a&[name=b] b
\ncline{a}{c}
\ncline{b}{c}
\end{psmatrix}
$$
\end{footnotesize}
\end{center}
\caption{The Hasse diagram for the poset $\Lambda$}
\label{lambda-fig}
\end{figure}

\subsection{The poset $\La$}

Can anything be said about the M\"obius function of $P^*$ if $P$ is
not a rooted forest?  Again, computer evidence suggests that the
answer is ``yes.''  Consider the poset $\La$ in
Figure~\ref{lambda-fig} which is the smallest one to which
Theorem~\ref{mufor} does not apply.  Let $T_n(x)$ denote the {\it
Tchebyshev polynomials of the first kind\/}, which can be 
defined as the unique polynomials such that
$$
T_n(\cos\theta)=\cos(n\theta).
$$

\bcon[Sagan-Vatter~\cite{sv:mfc}]
For all $i\le j$, $\mu(a^i,c^j)$ is the coefficient of $x^{j-i}$ in
$T_{i+j}(x)$.
\econ
Finding a proof of this conjecture by using generating functions or
any other means would be most welcome.

\bigskip

\noindent{\it Acknowledgment.}  We are indebted to Mihai Ciucu and
Andrew Sills for useful discussions about hypergeometric series.

\bigskip
\bibliographystyle{acm}
\begin{small}
\bibliography{ref}
\end{small}

\end{document}